\documentclass[12pt]{article}
\usepackage{amssymb,amsfonts,amsmath, psfrag,eepic,colordvi,graphicx,epsfig,ytableau}
\usepackage{amssymb,latexsym,graphics,array}
\usepackage{MnSymbol}
\usepackage[enableskew]{youngtab}
\usepackage{float,enumerate,enumitem}
\usepackage{tikz,ifthen}
\usepackage{pgfplots}
\usetikzlibrary{math}
\topskip  
\newcommand{\rmnum}[1]{\romannumeral #1}

\numberwithin{equation}{section}
\newtheorem{theorem}{Theorem}[section]

\newtheorem{conjecture}[theorem]{Conjecture}

\newtheorem{lemma}[theorem]{Lemma}

\begin{document}
	\parskip 6pt
	
	\pagenumbering{arabic}
	\def\sof{\hfill\rule{2mm}{2mm}}
	\def\ls{\leq}
	\def\gs{\geq}
	\def\SS{\mathcal S}
	\def\qq{{\bold q}}
	\def\MM{\mathcal M}
	\def\TT{\mathcal T}
	\def\EE{\mathcal E}
	\def\lsp{\mbox{lsp}}
	\def\rsp{\mbox{rsp}}
	\def\pf{\noindent {\it Proof.} }
	\def\mp{\mbox{pyramid}}
	\def\mb{\mbox{block}}
	\def\mc{\mbox{cross}}
	\def\qed{\hfill \rule{4pt}{7pt}}
	\def\block{\hfill \rule{5pt}{5pt}}
	\def\lr#1{\multicolumn{1}{|@{\hspace{.6ex}}c@{\hspace{.6ex}}|}{\raisebox{-.3ex}{$#1$}}}
	\def\red{\textcolor{red}}

	\begin{center}
		{\Large \bf On the  enumeration of  permutations avoiding chains of patterns}
	\end{center}
		\begin{center}
		{\small  Robin D.P. Zhou$^{*}$,\footnote{$^*$Corresponding author.}  \footnote{{\em E-mail address:} dapao2012@163.com.}  Yongchun Zang}
		
		College of Mathematics Physics and Information\\
		Shaoxing University\\
		Shaoxing 312000, P.R. China
	
	\end{center}
	
	\noindent {\bf Abstract.} 	
In 2019,  B\'ona and Smith introduced the notion of
strong pattern avoidance, saying that a permutation 
$\pi$ strongly avoids a pattern $\sigma$ if $\pi$ and $\pi^2$ both avoid $\sigma$.
Recently, Archer and Geary generalized the idea of strong pattern avoidance to chain avoidance, in which a permutation $\pi$ avoids a chain of patterns $(\tau^{(1)}:\tau^{(2)}:\cdots:\tau^{(k)})$ if the $i$-th power of the permutation  avoids the pattern $\tau^{(i)}$ for $1\leq i\leq k$.
In this paper, we give explicit formulae for the number of 
 sets of permutations avoiding certain chains of patterns.
Our results give affirmative answers to two conjectures  proposed by Archer and Geary.

	

	
	\section{Introduction}
	
	Let $\mathcal{S}_n$ denote the  symmetry group consisting of all permutations of $[n] = \{1,2,\ldots,n\}$, which we can view as bijections from $[n]$ to $[n]$.
A permutation $\pi$  can also be  viewed as a word	by associating  $\pi$  with the word $\pi(1)\pi(2)\cdots\pi(n)$.
	Throughout this paper, we always write a permutation  $\pi \in\mathcal{S}_n$ as a word
	$\pi = \pi_1\pi_2\cdots \pi_n$, where $\pi_i = \pi(i)$ for $1\leq i\leq n$.
	
Given a permutation $\pi \in \mathcal{S}_n$ and a permutation $\sigma \in \mathcal{S}_k$,
an {\em occurrence} of $\sigma$ in $\pi$ is a subsequence $\pi_{i_1}\pi_{i_2}\cdots \pi_{i_k}$
of $\pi$ that is order isomorphic to $\sigma$.
We say $\pi$   {\em contains}  the  pattern $\sigma$ if $\pi$ contains an occurrence of $\sigma$.
Otherwise, we say $\pi$ {\em avoids} the pattern $\sigma$  or
$\pi$ is {\em $\sigma$-avoiding}. 
For example, the permutation $1534627$ avoids the pattern $3142$ while it contains the  pattern $12345$ corresponding to the subsequence $13467$. 

The study of pattern avoidance can be traced back to the work of MacMahon
\cite{MacMahon} and has become a research focus in enumerative combinatorics over the past
half century.
A recent survey on permutation
patterns can be found in the books \cite{BonaBook} and \cite{Kitaev}.
The notion of pattern avoidance  has arisen many further variations.
B\'ona and Smith \cite{B.S} initiated the study of 
{\em strong pattern avoidance}, in which a permutation 
$\pi$ strongly {\em avoids} a pattern $\sigma$ if $\pi$ and $\pi^2$ both avoid $\sigma$.
The study of strong pattern avoidance has attracted the interest of combinatorial scholars such as \cite{B.C,JPan}.

Recently, Archer and Geary \cite{Archer1} generalized the idea of strong pattern avoidance to {\em chain avoidance}, in which a permutation $\pi$  {\em avoids} a chain (of patterns) $(\tau^{(1)}:\tau^{(2)}:\cdots:\tau^{(k)})$ if the $i$-th power 
of the  permutation (i.e., $\pi^i$) avoids the pattern $\tau^{(i)}$ for $1\leq i\leq k$.
Let $\mathcal{S}_n(\tau_{1}:\tau_{2}:\cdots:\tau_{k})$ denote the set of
$(\tau_{1}:\tau_{2}:\cdots:\tau_{k})$-avoiding permutations in $\mathcal{S}_n$. 
This definition can be  extended to include sets of patterns, which we will separate the patterns in a set with commas.
For example,  if $\pi$ avoids the chain 
$(\sigma, \rho:\tau)$, then 
$\pi$ avoids both $\sigma$ and $\rho$ and 
$\pi^2$ avoids $\tau$.
For example, let $\pi = 1325467$.
Then we have $\pi \in \mathcal{S}_7(231,1432:231)$ because
$\pi$ avoids both $231$ and $1432$, $\pi^2 = 1234567$ avoids $231$.

It would also be natural to consider {\em consecutive patterns} as part of the chain.
We say $\pi$ {\em contains} a consecutive pattern 
$\bar{\sigma}=\overline{\sigma_1\sigma_2\cdots\sigma_k}$ if there is an occurrence of $\sigma$ corresponding 
to  a consecutive subsequence $\pi_{i}\pi_{i+1}\cdots\pi_{i+k-1}$ ($1\leq i\leq n-k+1)$ of $\pi$.
Otherwise, we say $\pi$ {\em avoids} the (consecutive) pattern 
$\bar{\sigma}$ or $\pi$ is  {\em $\bar{\sigma}$-avoiding}.
For example, the permutation $\pi = 1534627$ avoids the pattern $\overline{321}$
while $\pi$ contains the pattern $\overline{213}$ corresponding to the consecutive subsequence $627$.

Archer and Geary \cite{Archer1} obtained the explicit formulae for the number of sets $\mathcal{S}_n(213,312:\tau)$ for $\tau\in\mathcal{S}_3$
 and further posed the following conjectures.

\begin{conjecture}\label{conj:sum231}
	For $n\geq 1$, we have $|\mathcal{S}_n(231,1432:231)| = L_{n+1}-\lceil\frac{n}{2}\rceil-1$, where $L_{n+1}$
	is the $(n + 1)$-th Lucas number (A000032, \cite{OEIS}).
\end{conjecture}

\begin{conjecture}\label{conj:sum213}
	For $n\geq 2$, we have $|\mathcal{S}_n(213,312:\overline{213})| = 2^{n-2}+n-1$. 
\end{conjecture}

The objective of this paper is to prove Conjecture \ref{conj:sum231} and 
Conjecture \ref{conj:sum213}.
Our proofs are based on the recursive decompositions of the permutations in $\mathcal{S}_n(231,1432:231)$ and the permutations in
$\mathcal{S}_n(213,312:\overline{213})$, respectively.


%
	
	\section{Avoiding the chain  $(231,1432:231)$}
	
%
	
This section is devoted to the proof of Conjecture \ref{conj:sum231}.
This will be accomplished by showing that the both sides of the equation in 
Conjecture \ref{conj:sum231} satisfy the same
recurrence and initial conditions.
The following  lemma is  crucial for our proof.

\begin{lemma}\label{lem:deco1}
	Let $n\geq 3$ and $\pi =\pi_1\pi_2\cdots \pi_n \in \mathcal{S}_n(231,1432:231)$.
	Then we have either $\pi_1 = n$, $\pi_n = n$ or $\pi_{n-1}\pi_n = n(n-1)$.
	
\end{lemma}

\pf 
The result is trivial for $n = 3$.
Now we let $n \geq 4$.
Assume that $\pi_1\neq n$, $\pi_n\neq n$ and $\pi_{n-1}\pi_n \neq n(n-1)$.
Let $\pi_k = n$. 
Then  $2\leq k \leq n-1$.
Combining the fact  $\pi$ is $231$-avoiding, we have that $\pi$ can be written in the form 
$\pi = \sigma n\tau$, where $\sigma$ and $\tau$ are not empty and all the elements in $\sigma$ are less than 
the elements in $\tau$.
Notice that $\pi_{n-1}\pi_n \neq n(n-1)$.
It is easily seen that $\tau$ contains at least two elements, namely, $k\leq n-2$.
We claim that the elements in  $\tau$ are increasing.
If not, there exists some $i>k$ such that $\pi_i>\pi_{i+1}$.
It follows that the subsequence $\pi_1n \pi_i \pi_{i+1}$ forms an 
occurrence of the pattern $1432$ in $\pi$, contradicting to the 
fact that $\pi$ is $1432$-avoiding.
This proves the claim.
Hence $\pi$ is of the form $\pi = \sigma nk(k+1)\cdots (n-1)$.
Now consider $\pi^2$. 
It is easily checked that $\pi^2(k) =\pi(\pi(k)) = \pi(n)=  n-1$, $\pi^2({k+1}) = \pi(\pi(k+1))=\pi(k) = n$ and 
$\pi^2(n) =\pi(\pi(n)) = \pi(n-1)= n-2$.
Then $\pi^2(k)\pi^2(k+1)\pi^2(n) = (n-1)n(n-2)$ forms an occurrence of the  pattern $231$ in $\pi^2$,
contradicting to the fact that $\pi^2$ is $231$-avoiding.
This completes the proof.
\qed

Given $n\geq 3$,
Lemma \ref{lem:deco1} enables us to divide the set $\mathcal{S}_n(231,1432:231)$
into the following three disjoint subsets:
\begin{align*}
	\mathcal{P}_n^1 &:= \{\pi \mid \pi=\pi_1\pi_2\cdots \pi_n \in \mathcal{S}_n(231,1432:231) \text{ and }  \pi_1 = n\},\\
    \mathcal{P}_n^2 &:= \{\pi \mid \pi=\pi_1\pi_2\cdots \pi_n \in \mathcal{S}_n(231,1432:231) \text{ and } \pi_n = n\},\\
    \mathcal{P}_n^3 &:= \{\pi \mid \pi=\pi_1\pi_2\cdots \pi_n \in \mathcal{S}_n(231,1432:231) \text{ and } \pi_{n-1}\pi_n = n(n-1)\}.
\end{align*}
We shall characterize the permutations in $\mathcal{P}_n^1$, $\mathcal{P}_n^2$ and 
$\mathcal{P}_n^3$, respectively.
The following theorem is needed for the characterization of $\mathcal{P}_n^1$.

\begin{theorem}\label{thm:yinyong}(\cite{B.S}, Theorem 3.1)
	For any permutation $\pi$ ending in $1$, the following two statements are equivalent.
\begin{itemize}
	\item [(A)] The permutation $\pi$ is strongly $312$-avoiding.
	\item [(B)] The permutation $\pi$ has form $\pi=(k+1)(k+2)\cdots nk(k-1)(k-2)\cdots1$ where $\lceil \frac{n}{2}\rceil \leq k\leq n-1$.
\end{itemize}
\end{theorem}	

\begin{lemma}\label{lem:P1}
	Let $n\geq 3$ and $\pi = \pi_1\pi_2\cdots \pi_n$.
	Then we have 
	\[\mathcal{P}_n^1 = \{\pi \mid \pi = n(n-1)\cdots (n-k+1)12\cdots (n-k),\, \lceil \frac{n}{2}\rceil \leq k\leq n-1\}.\]

\end{lemma}
\pf
Let $\mathcal{A}_n := \{\pi \mid \pi=\pi_1\pi_2\cdots \pi_n\in \mathcal{S}_n(312:312), \,\pi_n = 1\}$ and 
$\mathcal{B}_n :=\{\pi \mid \pi=\pi_1\pi_2\cdots \pi_n\in \mathcal{S}_n(231:231),\, \pi_1 = n\}$.
By the definition of $\mathcal{P}_n^1$, we have  $\mathcal{P}_n^1 \subseteq \mathcal{B}_n$.
As the pattern $231$ is the inverse of $312$, 
it is straightforward to see that $\pi$ is strongly $231$-avoiding
if and only if $\pi^{-1}$ is strongly $312$-avoiding.
Note that a permutation $\pi$  begins with  $n$ if and only if $\pi^{-1}$ ends in  $1$.
Therefore, the inverse map is a bijection between 
$\mathcal{A}_n$ and $\mathcal{B}_n$.
From Theorem \ref{thm:yinyong}, we obtain that 
\[\mathcal{A}_n = \{\pi \mid \pi = (k+1)(k+2)\cdots nk(k-1)(k-2)\cdots1,\, \lceil \frac{n}{2}\rceil \leq k\leq n-1\}.\]
It can be checked that the permutation $(k+1)(k+2)\cdots nk(k-1)(k-2)\cdots1$ is the inverse of $n(n-1)\cdots (n-k+1)12\cdots (n-k)$.
It yields that
	\[\mathcal{B}_n = \{\pi \mid \pi = n(n-1)\cdots (n-k+1)12\cdots (n-k),\, \lceil \frac{n}{2}\rceil \leq k\leq n-1\}.\]
	Observe that all the permutations of the form $n(n-1)\cdots (n-k+1)12\cdots (n-k)$ are $1432$-avoiding.
	We derive that $\mathcal{B}_n\subseteq \mathcal{P}_n^1$.
	Hence we have $\mathcal{P}_n^1 = \mathcal{B}_n$.
	This completes the proof.
	\qed
	
		Given $\sigma\in \mathcal{S}_{k}$ and $\tau\in \mathcal{S}_{m}$, let  $\sigma\oplus\tau$ denote the  {\em direct sum} of $\sigma$ and $\tau$ defined by 
	
	\[
	(\sigma\oplus\tau)(i) = 
	\begin{cases}
		\sigma(i), \text{\quad if $1\leq i\leq k$}, \\
		\tau(i-k)+k,  \text{\quad if $k+1\leq i\leq k+m$.}
	\end{cases}
	\]
	For example, 231 $\oplus$ 21 = 23154. 
	The following fact is straightforward from the definition 
	of direct sum.
	\begin{align}\label{fact:1}
	(\sigma\oplus \tau)^2 = \sigma^2 \oplus \tau^2.
	\end{align}


%
%

\begin{lemma}\label{lem:P2}
	Let $n\geq 3$ and $\pi = \pi_1\pi_2\cdots \pi_n\in \mathcal{S}_n$.
	Then $\pi \in \mathcal{P}_n^2$ if and only if 
	$\pi = \sigma\oplus1$ for some  $\sigma\in\mathcal{S}_{n-1}(231,1432:231)$.
\end{lemma}

\pf
If  $\pi\in\mathcal{P}_n^2$, then we have $\pi_n = n$.
By the definition of the direct sum, we have that $\pi = \sigma \oplus 1$ for some $\sigma \in \mathcal{S}_{n-1}$.
From (\ref{fact:1}) we have $\pi^2 = \sigma^2 \oplus 1$.
This means that $\sigma$ (resp. $\sigma^2$) is a subsequence of
$\pi$ (resp. $\pi^2$).
Since $\pi$ avoids the chain $(231,1432:231)$, 
we deduce that $\sigma$ also avoids the chain $(231,1432:231)$,
that is, $\sigma\in\mathcal{S}_{n-1}(231,1432:231)$.

Conversely,  if $\pi = \sigma\oplus1$ for some  $\sigma\in\mathcal{S}_{n-1}(231,1432:231)$, 
we have $\pi_n = n$ and $\pi^2 = \sigma^2 \oplus 1$.
To prove $\pi \in \mathcal{P}_n^2$,
it suffices to  show that $\pi$ avoids the chain $(231,1432:231)$.
Observe that the element $\pi_n = n$ can not appear in any occurrences 
of the pattern $231$ (resp. $1432$) in $\pi$.
Combining the fact that $\sigma$ avoids the chain $(231,1432:231)$,
we deduce that $\pi$ avoids both $231$ and $1432$.
By a similar argument as above, we have  $\pi^2$ is $231$-avoiding.
In conclusion, we have $\pi$ avoids the chain $(231,1432:231)$,
as desired.
\qed

\begin{lemma}\label{lem:P3}
	Let $n\geq 3$ and $\pi = \pi_1\pi_2\cdots \pi_n\in \mathcal{S}_n$.
	Then $\pi \in \mathcal{P}_n^3$ if and only if 
	$\pi = \sigma\oplus21$ for some  $\sigma\in\mathcal{S}_{n-2}(231,1432:231)$.
\end{lemma}
	
	\pf
	If  $\pi\in\mathcal{P}_n^3$, then we have $\pi_{n-1}\pi_n = n(n-1)$.
	By the definition of the direct sum, we have that $\pi = \sigma \oplus 21$ for some $\sigma \in \mathcal{S}_{n-2}$.
	Then we obtain $\pi^2 = \sigma^2 \oplus 12$ from  (\ref{fact:1}).
	This means that $\sigma$ (resp. $\sigma^2$) is a subsequence of
	$\pi$ (resp. $\pi^2$).
	Since $\pi$ avoids the chain $(231,1432:231)$, 
	we deduce that $\sigma$ also avoids the chain $(231,1432:231)$,
	that is, $\sigma\in\mathcal{S}_{n-2}(231,1432:231)$.
	
	Conversely,  assume that $\pi = \sigma\oplus 21$ for some  $\sigma\in\mathcal{S}_{n-2}(231,1432:231)$.
	By the definition of direct sum and (\ref{fact:1}),
	we have $\pi_{n-1}\pi_n = n(n-1)$ and $\pi^2 = \sigma^2 \oplus 12$.
	To prove $\pi \in \mathcal{P}_n^3$,
	it suffices to  show that $\pi$ avoids the chain $(231,1432:231)$.
	As $\pi_{n-1}\pi_n = n(n-1)$, 
	we deduce that neither   $n$ nor $n-1$ can 
	 appear in any occurrences 
	of the pattern $231$ (resp. $1432$) in $\pi$.
	Combining the fact that $\sigma$ avoids the chain $(231,1432:231)$,
	we deduce that $\pi$ avoids both $231$ and $1432$.
	By a similar argument as above, we have  $\pi^2$ is $231$-avoiding.
	Therefore  $\pi$ avoids the chain $(231,1432:231)$,
	as desired.
	\qed

\noindent  {\bf Proof of Conjecture \ref{conj:sum231}.} 	
Let $n \geq 1$ and $f(n) = |\mathcal{S}_n(231,1432:231)|$.
Then we have 
\[f(n) = |\mathcal{P}_n^1|+|\mathcal{P}_n^2|+|\mathcal{P}_n^3|\]
for $n\geq 3$.
Lemmas \ref{lem:P1}, \ref{lem:P2} and \ref{lem:P3} tell us that 
$|\mathcal{P}_n^1| = n-\lceil \frac{n}{2}\rceil =\lceil\frac{n-1}{2}\rceil$,
$|\mathcal{P}_n^2| = f(n-1)$ and $|\mathcal{P}_n^3| = f(n-2)$,  respectively.
Therefore we have 
\[f(n)=\lceil\frac{n-1}{2}\rceil+f(n-1)+f(n-2)\] 
for $n\geq 3$ with the initial conditions $f(1) = 1$ and $f(2) = 2$.
Notice that the Lucas numbers $L_n$ satisfy that 
$L_n = L_{n-1} +L_{n-2}$ for $n\geq 3$ with the initial conditions
$L_1 = 1$ and $L_2 = 3$.
It is routine to check that $f(n)$ and $L_{n+1}-\lceil\frac{n}{2}\rceil-1$ satisfy the 
same recurrence and initial conditions.
Thus 
$f(n) = L_{n+1}-\lceil\frac{n}{2}\rceil-1$ for $n\geq 1$.
This completes the proof.
\qed

\section{Avoiding the chain $(213,312:\overline{213})$}	

This section is devoted to the proof of Conjecture 
\ref{conj:sum213}.
To this end, we need to give a characterization of the permutations in $\mathcal{S}_n(213,312:\overline{213})$.
Note that the permutations in $\mathcal{S}_n(213,312)$
are exactly the set of unimodal permutations,
i.e.,  
those permutations $\pi = \pi_1\pi_2\cdots \pi_n$  with $\pi_1<\pi_2<\cdots <\pi_k>\pi_{k+1}>\cdots > \pi_n$ for some $1\leq k\leq n$ \cite{Gannon}.
For  a unimodal permutation $\pi = \pi_1\pi_2\cdots \pi_n$, we have either $\pi_1 = 1$ or 
$\pi_n=1$.
Hence, for $n\geq 2$, we can divide the set  $\mathcal{S}_n(213,312:\overline{213})$  into
the following two disjoint subsets:
\begin{align*}
	\mathcal{Q}_n^1 &:= \{\pi \mid \pi=\pi_1\pi_2\cdots \pi_n \in \mathcal{S}_n(213,312:\overline{213}) \text{ and }  \pi_1 = 1\},\\
	\mathcal{Q}_n^2 &:= \{\pi \mid \pi=\pi_1\pi_2\cdots \pi_n \in \mathcal{S}_n(213,312:\overline{213}) \text{ and } \pi_n = 1\}.
\end{align*}

We will give characterizations of the unimodal permutations in $\mathcal{Q}_n^1$ and $\mathcal{Q}_n^2$, respectively.

\begin{lemma}\label{lem:start1}
	Let $n\geq 3$ and $\pi =\pi_1\pi_2\cdots \pi_n\in \mathcal{S}_n$.
	Then we have $\pi \in \mathcal{Q}_n^1$ if and only if 
	$\pi=1\oplus\sigma$ for some $\sigma\in\mathcal{S}_{n-1}(213,312:\overline{213})$.  
\end{lemma}

\pf
If $\pi \in \mathcal{Q}_n^1$, then we have $\pi_1 = 1$.
By the definition of direct sum, 
the permutation $\pi$ can be written as $\pi = 1\oplus \sigma$ for 
some $\sigma\in \mathcal{S}_{n-1}$.
From (\ref{fact:1}), we have $\pi^2 = 1\oplus \sigma^2$.
Hence $\sigma$ (resp. $\sigma^2$) is order isomorphic to 
$\pi_2\pi_3\cdots \pi_n$ (resp. $\pi^2(2)\pi^2(3)\cdots \pi^2(n)$).
Then $\sigma\in\mathcal{S}_{n-1}(213,312:\overline{213})$ follows directly from the fact that $\pi$ avoids the chain 
$(213,312:\overline{213})$.

Conversely, if 
$\pi=1\oplus\sigma$ for some $\sigma\in\mathcal{S}_{n-1}(213,312:\overline{213})$.
Then we have $\pi_1= 1$ and $\pi^2 = 1\oplus \sigma^2$.
To prove $\pi \in \mathcal{Q}_n^1$,
it suffices to  show that $\pi$ avoids the chain $(213,312:\overline{213})$.
Observe that the element $\pi_1=1$ can not appear in any occurrences 
of the pattern $213$ (resp. $312$) in $\pi$.
Combining the fact that $\sigma$ avoids the chain $(213,312:\overline{213})$,
we deduce that $\pi$ avoids both $213$ and $312$.
By a similar argument as above, we have  $\pi^2$ is $\overline{213}$-avoiding.
Therefore the permutation $\pi$ avoids the chain $(213,312:\overline{213})$,
as desired.
\qed

The following lemma is needed for the characterization of 
permutations in $\mathcal{Q}_n^2$.

\begin{lemma}\label{lem:peak}
	Let $n\geq 3$ and  $\pi = \pi_1\pi_2\cdots \pi_n\in\mathcal{S}_n(213,312)$.
	If $\pi^2(i-1)\pi^2(i)\pi^2(i+1)$  is an occurrence of 
	the consecutive pattern $\overline{213}$ in $\pi^2$, then we have
	$\pi_i = n$. 
\end{lemma}

\pf 
We have known that the permutations in $\mathcal{S}_n(213,312)$ are unimodal permutations.
If $\pi_i \neq n$,
 then we have either $\pi_{i-1}<\pi_i<\pi_{i+1}$ or $\pi_{i-1}>\pi_i>\pi_{i+1}$.
In the former case,  we have $\pi(\pi_{i-1})\pi(\pi_i)\pi(\pi_{i+1}) = \pi^2(i-1)\pi^2(i)\pi^2(i+1)$
is a subsequence of $\pi$.
Since $\pi^2(i-1)\pi^2(i)\pi^2(i+1)$ is an occurrence of 	the consecutive pattern $\overline{213}$ in $\pi^2$, 
we obtain that $\pi$ contains the pattern $213$, 
contradicting to the fact that $\pi$ is $213$-avoiding.
In the later case, we have $\pi(\pi_{i+1})\pi(\pi_i)\pi(\pi_{i-1}) = \pi^2(i+1)\pi^2(i)\pi^2(i-1)$
is a subsequence of $\pi$.
Since $\pi^2(i-1)\pi^2(i)\pi^2(i+1)$ is an occurrence of 	the consecutive pattern $\overline{213}$ in $\pi^2$, 
we have $\pi^2(i+1)\pi^2(i)\pi^2(i-1)$ is an occurrence of 
the pattern $312$ in $\pi$, 
contradicting to the fact that $\pi$ is $312$-avoiding.
This completes the proof.
\qed

\begin{lemma}\label{lem:end1}
For $n\geq 3$, we have 
	\[\mathcal{Q}_n^2 =\{\pi = \pi_1\pi_2\cdots \pi_n \mid 
	\pi =23\cdots n1  \text{ or } \pi = \sigma n(n-1)\tau 1\},
	\]
	where $\sigma$ (possibly empty) is increasing and $\tau$ (possibly empty) is decreasing.
\end{lemma}

\pf 
Let $\pi = \pi_1\pi_2\cdots \pi_n$ be a unimodal permutation with $\pi_k = n$ and $\pi_n=1$.
Then  we have either 
$\pi_{k-1} = n-1$ or $\pi_{k+1} = n-1$.

\noindent {\bf Case (\rmnum{1}).} 
$\pi_{k-1}=n-1$. \\
First  we consider the permutation $\pi=23\cdots n1$. 
We have $\pi^2=34\cdots n12$. 
By the definition of $\mathcal{Q}_n^2$, one can easily check that $\pi\in\mathcal{Q}_n^2$. 
Now assume that $\pi \neq 23\cdots n1$.
It yields that there exists at least one element between 
$n$ and $1$ in $\pi$.
We claim that $\pi \nin \mathcal{Q}_n^2$.
We will prove the claim by showing that 
$\pi^2(k-1)\pi^2(k)\pi^2(k+1)$  is an occurrence of 
the consecutive pattern $\overline{213}$ in $\pi^2$.
Since $\pi$ is unimodal with $\pi_k =n$ and $\pi_n = 1$, we have $\pi_i>i $ for $i\in[1,k]$ and $\pi_{k+1}\ge\pi_j\ge\pi_{n-1}$ for $j\in[k+1,n-1]$. 
It is routine to check that $\pi^2(k) = \pi(\pi_k) = \pi_n=1$ and $\pi^2(k-1) = \pi(\pi_{k-1}) = \pi_{n-1}$.
If $\pi_{k+1}\leq k$, we have $\pi^2(k+1)=\pi(\pi_{k+1})>\pi_{k+1}$.
As $\pi_{k+1} \geq \pi_{n-1} = \pi^2(k-1)$, 
we have $\pi^2(k+1) > \pi^2(k-1)$.
It follows that  
$\pi^2(k-1)\pi^2(k)\pi^2(k+1)$ is an occurrence of the pattern  $\overline{213}$ in $\pi^2$,  as desired.
Now suppose that $k < \pi_{k+1} < n-1$. 
We have $\pi^2(k-1)=\pi(n-1)<\pi(\pi_{k+1})=\pi^2(k+1)$.
 Again $\pi_{k-1}^2\pi_k^2\pi_{k+1}^2$ is an occurrence of the pattern  $\overline{213}$ in $\pi^2$, as desired.
 Hence, the claim is verified.
 In conclusion, there exists only one permutation 
 $\pi=23\cdots n1\in \mathcal{Q}_n^2$ if $\pi_{k-1} = n-1$.

\noindent {\bf Case (\rmnum{2}).} 
$\pi_{k+1}=n-1$. \\
Recall that $\pi$ is unimodal.
 In this case, the permutation $\pi$ can be written in the form $\pi=\sigma n(n-1)\tau1$ where $\sigma$ (possibly empty) is increasing and $\tau$ (possibly empty) is decreasing.
 We need to show that all such permutations are contained in the set $\mathcal{Q}_n^2$.
 If $\sigma$ is empty,  it is straightforward to check that 
 $\pi = n(n-1)\cdots 1 \in \mathcal{Q}_n^2$, as desired.
 Now assume that $\sigma$ is nonempty.  
 To prove $\pi \in \mathcal{Q}^2_n$, we need to show that 
 $\pi^2$ avoids the consecutive pattern $\overline{213}$.
 From Lemma \ref{lem:peak},  it is sufficient to 
 show that $\pi^2(k-1)\pi^2(k)\pi^2(k+1)$  is not an occurrence of  $\overline{213}$ in $\pi^2$.
Observe that  $k\leq\pi_{k-1} <n-1$.
Since $\pi$ is unimodal, we have  $\pi^2(k-1)=\pi(\pi_{k-1})>\pi(n-1)=\pi^2(k+1)$.
Thus, $\pi^2(k-1)\pi^2(k)\pi^2(k+1)$ is not an occurrence 
of $\overline{213}$ in $\pi^2$, as desired.
\qed

	\noindent  {\bf Proof of Conjecture \ref{conj:sum213}.} 	
	Let $n \geq 2$ and $g(n) = |\mathcal{S}_n(213,312:\overline{213})|$.
	Then we have 
	\[g(n) = |\mathcal{Q}_n^1|+|\mathcal{Q}_n^2|.\]
	By Lemma \ref{lem:start1}, we have $|\mathcal{Q}_n^1| = g(n-1)$ for $n\geq 3$.
Lemma \ref{lem:end1} tells us that 	\[\mathcal{Q}_n^2 =\{\pi = \pi_1\pi_2\cdots \pi_n \mid 
\pi =23\cdots n1  \text{ or } \pi = \sigma n(n-1)\tau 1\},
\]
where $\sigma$ (possibly empty) is increasing and $\tau$ (possibly empty) is decreasing.
A permutation $\pi\in \mathcal{Q}_n^2$ with the form  $\pi= \sigma n(n-1)\tau 1$ is uniquely determined by 
the set consisting of the elements in $\sigma$.
We can choose the elements of $\sigma$ in $\binom{n-3}r$ way if $\sigma$ contains $r$ elements.
Then  we deduce that
\[|\mathcal{Q}_n^2| = 1+ \sum_{\mathrm{r}=0}^{n-3}\binom{n-3}r=2^{n-3}+1\]
for $n\geq 3$.
	Therefore we have 
	\[g(n)=g(n-1)+2^{n-3}+1\] 
	for $n\geq 3$.
	Taking into account the initial conditions $g(2)= 2$, the unique solution is $g(n) =  2^{n-2}+n-1$.
	This completes the proof.
	\qed

	\section*{Acknowledgments}
	The work  was supported by
	the National Natural
	Science Foundation of China (11801378).
	

\end{document}